\documentclass[11pt]{amsart}

\usepackage{times,epsf,amssymb,amsmath,rlepsf}
 \usepackage[ps2pdf, bookmarks, colorlinks]{hyperref}
\hypersetup{
    colorlinks,%
    citecolor=blue,%
    linkcolor=red,%
}

\begin{document}

\newtheorem{thm}{Theorem}[section]
\newtheorem{lem}[thm]{Lemma}
\newtheorem{cor}[thm]{Corollary}

\theoremstyle{definition}
\newtheorem{defn}{Definition}[section]

\theoremstyle{remark}
\newtheorem{rmk}{Remark}[section]

\def\square{\hfill${\vcenter{\vbox{\hrule height.4pt \hbox{\vrule
width.4pt height7pt \kern7pt \vrule width.4pt} \hrule height.4pt}}}$}

\def\T{\mathcal T}

\newenvironment{pf}{{\it Proof:}\quad}{\square \vskip 12pt}

\title{Asymptotic Plateau Problem}
\author{Baris Coskunuzer}
\address{Department of Mathematics \\ Koc University \\ Sariyer, Istanbul 34450 Turkey}
\email{bcoskunuzer@ku.edu.tr}
\thanks{The author is partially supported by EU-FP7 Grant IRG-226062, TUBITAK Grant 109T685 and TUBA-GEBIP Award.}

\maketitle


\newcommand{\Si}{S^2_{\infty}(\mathbf{H}^3)}
\newcommand{\PI}{\partial_{\infty}}
\newcommand{\SI}{S^{n}_{\infty}(\mathbf{H}^{n+1})}
\newcommand{\BHH}{\mathbf{H}^{n+1}}
\newcommand{\CH}{\mathcal{C}(\Gamma)}
\newcommand{\BH}{\mathbf{H}}
\newcommand{\BR}{\mathbf{R}}
\newcommand{\BC}{\mathbf{C}}
\newcommand{\BZ}{\mathbf{Z}}

\begin{abstract}

This is a survey of old and recent results about the asymptotic Plateau problem. Our aim is to give a fairly complete
picture of the field, and present the current situation.

\end{abstract}

\section{Introduction}

The asymptotic Plateau problem in hyperbolic space basically asks the existence of an area
minimizing submanifold $\Sigma \subset \BHH$ asymptotic to given submanifold $\Gamma\subset \SI$.
In this survey article, we will cover old and recent results on the problem. Most of the time, we
will give the essential ideas of the proofs. Our aim is to give a nice expository introduction for
the interested researchers, and to present a picture of this growing field.

\tableofcontents

\section{Preliminaries}

In this section, we will overview the basic results which we will use in the following sections.
First, we will give the definitions of area minimizing surfaces. First set of the definitions are
about compact submanifolds. The second set of the definitions are their generalizations to the
noncompact submanifolds.

\begin{defn}(Compact Case) Let $D$ be a compact disk in a manifold $X$. Then, $D$ is an {\em area minimizing disk} in $X$ if
$D$ has the smallest area among the disks in $X$ with the same boundary. Let $S$ be a compact
submanifold with boundary in a manifold $X$. Then, $S$ is an {\em absolutely area minimizing
submanifold} in $X$ if $S$ has the smallest volume among all submanifolds (no topological
restriction) with the same boundary in $X$. The absolutely area minimizing surfaces and
hypersurfaces can be defined likewise.
\end{defn}

\begin{defn}(Noncompact Case) An {\em area minimizing plane} (least area plane) is a complete plane in a manifold $X$
such that any compact subdisk in the plane is an area minimizing disk in $X$. Let $\Delta$ be a
complete submanifold in a manifold $X$. Then, $\Delta$ is an {\em absolutely area minimizing
submanifold} in $X$ if any compact part (codimension-$0$ submanifold with boundary) of the $\Delta$
is an absolutely area minimizing hypersurface in $X$. The absolutely area minimizing surfaces,
hypersurfaces and hyperplanes can be defined likewise.
\end{defn}

\begin{defn} A {\em minimal surface (submanifold or hypersurface)} in a manifold $X$ is a surface
(submanifold or hypersurface) whose mean curvature vanishes everywhere.
\end{defn}

Note that the mean curvature being $0$ is equivalent to be locally area minimizing \cite{CM1}.
Hence, all area minimizing surfaces and hypersurfaces are also minimal.

\begin{defn} (Convex Hull) Let $A$ be a subset of $\SI$. Then the \textit{convex hull} of $A$, $CH(A)$,
is the smallest closed convex subset of $\BHH$ which is asymptotic to $A$. Equivalently, $CH(A)$
can be defined as the intersection of all supporting closed half-spaces of $\BHH$ \cite{EM}.
\end{defn}

Note that $\PI(CH(A)) = A$ for any $A\subset \SI$ (Note that this is a special property of $\BHH$,
see \cite{HLS}). In general, we say a subset $\Sigma$ of $X$ has the convex hull property if it is
in the convex hull of its boundary in $X$, i.e. $\Sigma \subset CH(\partial\Sigma)$. In special
case, if $\Sigma$ is a complete and noncompact hypersurface in $\BHH$, then we say $\Sigma$ has
convex hull property if it is in the convex hull of its asymptotic boundary, i.e. $\Sigma\subset
CH(\PI\Sigma)$. The minimal hypersurfaces in $\BHH$ have convex hull property.

\begin{lem} \cite{An1}
Let $\Sigma$ be a minimal submanifold in $\BHH$ with $\PI \Sigma = \Gamma$. Then $\Sigma \subset
CH(\Gamma)$.
\end{lem}

The idea is quite simple. Let $\Sigma$ be a minimal submanifold in $\BHH$ with $\PI\Sigma = \Gamma$. Let $K$ be a
nonsupporting halfspace in $\BHH$, i.e. $\PI K \cap \Gamma = \emptyset$. Since $K$ is a halfspace in $\BHH$, we can
foliate $K$ with geodesic planes whose asymptotic boundaries are in $\PI K$. Then, by maximum principle \cite{CM1}, $K
\cap \Sigma = \emptyset$, and hence $\Sigma \subset CH(\Gamma)$. We should also note that instead of smooth
submanifolds, if one deals with area minimizing rectifiable currents, or stationary varifolds, which might have some
singularities, for this type of results, one needs {\em strong maximum principle} results which applies to these
settings due to Simon \cite{Si}, Solomon-White \cite{SW}, Ilmanen \cite{Il} and Wickramasekera \cite{Wi}.

Throughout the paper, $\BHH$ will represent the hyperbolic $n+1$-space. $\BHH$ has a natural
compactification $\overline{\BHH} = \BHH \cup \SI$ where $\SI$ is the sphere at infinity of $\BHH$.
If $\Sigma$ is a subset of $\BHH$, the asymptotic boundary of $\Sigma$, say $\PI \Sigma$, can be
defined as $\PI \Sigma = \overline{\Sigma}-\Sigma$ where $\overline{\Sigma}$ is closure of $\Sigma$
in $\overline{\BHH}$ in the Euclidean metric. In the remaining of the paper, it is mostly a good
idea to imagine $\BHH$ in the Poincare ball model.

\section{Existence}

There are basically $2$ types of existence results for the asymptotic Plateau problem. The first
type is the existence of absolutely area minimizing submanifolds in $X$ for a given asymptotic
boundary in $\PI X$. In this type, there is no topological restriction on the submanifolds. The
other type is the fixed topological type. The area minimizing submanifold with the given asymptotic
data should also be in the given topological type.

\subsection{Absolutely Area Minimizing Submanifolds}\ \\

By using geometric measure theory methods, Michael Anderson solved the asymptotic Plateau problem for absolutely area
minimizing varieties for any dimension and codimension in \cite{An1}.

\begin{thm} \cite{An1} Let $\Gamma^p \to \SI$ be an embedded closed submanifold in the sphere at infinity of $\BHH$. Then there exists a
complete, absolutely area minimizing locally integral p+1-current $\Sigma$ in $\BHH$ asymptotic to $\Gamma^p$ at
infinity.
\end{thm}

\begin{pf} (Sketch) Let $\Gamma^p$ be an embedded closed submanifold in $\SI$. First, Anderson proves a monotonicity formula for
stationary p+1-currents such that the ratio between the volume of a stationary p+1-current restricted to a $r$-ball in
$\BHH$ and the volume of p+1-dimensional $r$-ball is nondecreasing in $r$ (\cite{An1}, Theorem 1). Then, he defines a
sequence of closed submanifolds $\Gamma_t^p$ in $\BHH$ such that $\Gamma_t^p \subset \partial B_t(0)$ and $\Gamma_t^p
\to \Gamma^p$.

Let $\Sigma_t$ be an area minimizing integral $p$-current with $\partial \Sigma_t = \Gamma_t$ \cite{Fe}. Then by using
the monotonicity formula, he gives a lower bound for the volume of $\Sigma_t$ restricted to $r$-ball, i.e. $c_r< ||
\Sigma_t|_{B_r} ||$. Also, by using the area minimizing property of $\Sigma_t$, he easily gives an upper bound $C_r$
for the volume of $\Sigma_t$ restricted to $r$-ball. Then, $c_r< || \Sigma_t|_{B_r} || < C_r$. Hence, by using
compactness theorem for integral currents (See \cite{Fe}, \cite{Mo}), he gets a convergent subsequence for
$\{\Sigma_i\}$ \textit{for each} $r$-ball. Then, by using diagonal subsequence argument, he extracts a convergent
subsequence $\Sigma_{i_j} \to \Sigma$ where $\Sigma$ is an area minimizing integral p+1-current with $\PI \Sigma =
\Gamma ^p$.
\end{pf}

\begin{rmk} This result is one of the most important results in the field. This seminal paper can be considered as the beginning of
the study of the asymptotic Plateau problem. Later, we will see various generalizations of this result to different
settings. Note that the embeddedness assumption on the given asymptotic boundary is very essential. In \cite{La2}, Lang
constructed immersed examples in $\SI$ with no solutions to asymptotic Plateau problem.
\end{rmk}

\begin{rmk} (Interior Regularity) By interior regularity results of geometric measure theory \cite{Fe}, \cite{Mo}, when $p=n-1$, the currents in
theorem are smoothly embedded hypersurfaces except for a singular set of Hausdorff dimension $n-7$. In particular when
$p=n-1<6$, $\Sigma$ is a smoothly embedded hypersurface in $\BHH$. In the higher codimension case ($p<n-1$), the
interior regularity results say that the absolutely area minimizing currents are smoothly embedded p+1-submanifolds in
$\BHH$ except for a singular set of Hausdorff dimension $p-1$.
\end{rmk}

Later, again by using geometric measure theory methods, Lang and Bangert generalized this result to Gromov hyperbolic
Hadamard manifolds with bounded geometry, and some other special cases in \cite{La1}, \cite{La2}, and \cite{BL} (See
also \cite{Gr}).

\begin{thm} \cite{La1}
Let $X$ be a Gromov hyperbolic Hadamard $n$-manifold with bounded geometry. Let $\Gamma$ be a $p$ dimensional closed
submanifold in $\PI X$. Then there exists a complete, absolutely area minimizing locally integral p+1-current $\Sigma$
in $X$ asymptotic to $\Gamma^p$ at infinity.
\end{thm}

Note that the varieties constructed in theorems above are absolutely area minimizing, and has no topological
restrictions on them. Another interesting case is the fixed topological type.

\subsection{Fixed Topological Type}\ \\

In above result, Anderson got absolutely area minimizing varieties asymptotic to given submanifold
in the asymptotic sphere. As there is no topological restriction on the objects, we have no idea
about the topological properties of them.

Another interesting case in Plateau problem is the fixed topological type. The question is to find
the smallest area surface in the given topological type with the given boundary. Its generalization
to the asymptotic Plateau problem is natural.

On the other hand, hyperbolic $3$-manifolds, and essential $2$-dimensional submanifolds in them are
very active research area. By essential, we mean $\pi_1$-injective surfaces, and they are very
important tools to understand the structure of the hyperbolic manifold by using geometric topology
tools. At this point, when we pass to the universal cover of the hyperbolic manifold and essential
surfaces in them, the asymptotic Plateau problem in disk type becomes an important technique for
construction of area minimizing representative of these essential surfaces in $3$-manifolds.

In \cite{An2}, Anderson focused on the asymptotic Plateau problem in disk type, and gave an
existence result in dimension $3$.

\begin{thm} \cite{An2} Let $\Gamma$ be a simple closed curve in $\Si$. Then, there exists a complete, area minimizing plane
$\Sigma$ in $\BH^3$ with $\PI \Sigma=\Gamma$.
\end{thm}

The proof is very similar to the proof of the previous theorem.  The basic difference is instead of
using area minimizing surfaces $\{\Sigma_t\}$ with $\partial \Sigma_t = \Gamma_t$, he used the area
minimizing {\em embedded} disks $\{D_t\}$ with $\partial D_t =\Gamma_t$. The existence of the disks
comes from the solution of Plateau problem in disk type. However, the essential point is that the
disks are embedded and they are given by \cite{AS} (later more general result given by \cite{MY1}).
Hence, by using similar ideas, Anderson extracted a limit $D_i \to \Sigma$ where $\Sigma$ is an
area minimizing plane in $\BHH$ with $\PI \Sigma = \Gamma$.

\begin{rmk} Note that this result is for just dimension $3$, it is not known if its generalization to higher dimensions is true or not.
It might be possible to construct area minimizing hyperplanes in $\BHH$ for any dimension, by generalizing these ideas
and White's results in \cite{Wh} to replace the sequence of disks $\{D_i\}$ in Anderson's proof with compact area
minimizing hyperplanes in $\BHH$.
\end{rmk}

Also, in \cite{An2}, Anderson constructed special Jordan curves in $\Si$ such that the absolutely
area minimizing surface given by Theorem 3.1 cannot be a plane (\cite{An2}, Theorem 4.5). Indeed,
he constructed examples with genus $g>g_0$ for any given genus $g_0$. He also used these surfaces
for some nonuniqueness results which we mention later.

In the same context, de Oliveira and Soret showed the existence of a complete stable minimal
surface in $\BH^3$ for any given topological type of a surface with boundary. Also, they studied
the isotopy type of these surfaces in some special cases. The main difference with Anderson's
existence result is that Anderson starts with the asymptotic boundary data, and gives an area
minimizing hypersurface where there is no control on the topological type, while de Oliveira and
Soret starts with a surface with boundary and constructs a stable minimal embedded surface of this
type whose asymptotic boundary is essentially determined by the surface.

\begin{thm} \cite{OS}
Let $M$ be a compact orientable surface with boundary. Then $int(M)$ can be minimally, completely,
properly and stably embedded in $\BH^3$. Furthermore, the embeddings extends smoothly to an
embedding from $M$ to $\overline{\BH^3}$, the compactification of $\BH^3$.
\end{thm}

On the other hand, Gabai gave another construction for Theorem 3.3. Indeed, he needed this results
for more general metrics, and he gave a topological construction for such area minimizing planes in
more general settings.

\begin{thm} \cite{Ga}
Let $X$ be $\BH^3$ with a different Riemannian metric induced from a metric on a closed
$3$-manifold. Let $\Gamma$ be a simple closed curve in $S^2_\infty(X)$. Then, there exists a
$D^2$-limit lamination $\sigma$ whose leaves are area minimizing planes in $X$ with $\PI \sigma
=\Gamma$.
\end{thm}

\begin{pf} (Sketch) Let $X$ be $\widetilde{M}$ where $\widetilde{M}$ is the universal cover of a hyperbolic $3$-manifold $M$ with any
Riemannian metric. In a similar fashion to the Anderson's proof, Gabai starts with a sequence of area minimizing disks
$\{D_i\}$ in $X$ with $\partial D_i = \Gamma_i \to \Gamma$. To get a limiting plane here, instead of using the
compactness theorem of geometric measure theory, he extracts some kind of Gromov-Hausdorff limit $\sigma$ of the
sequence $\{D_i\}$ by using minimal surface tools and techniques of \cite{HS}. In particular, the sequence $\{D_{i}\}$
of embedded disks in a Riemannian manifold $X$ \textit{converges} to the lamination $\sigma$ if

i) For any convergent sequence $\{x_{n_i} \}$ in $X$ with  $x_{n_i} \in D_{n_i}$ where $n_i$ is a strictly increasing
sequence, $\lim x_{n_i} \in \sigma$.

ii) For any $x \in \sigma$, there exists a sequence $\{x_i\}$ with $x_i \in D_i$ and $\lim x_i = x$ such that there
exist embeddings $f_{i}: D^2 \to D_i$ which converge in the $C^{\infty }$-topology to a smooth embedding $f:D^2 \to
L_{x}$, where $x_i \in f_{i}(Int(D^2))$, and $L_{x}$ is the leaf of $\sigma $ through $x$, and $x\in f(Int(D^2))$.

We call such a lamination $\sigma$ a \textit{$D^2$-limit lamination}. Here, the \textit{topological} limit $\sigma$ is
essentially all the limit points of a very special subsequence. Then, since locally these are limits of area minimizing
disks, by using the techniques of \cite{HS} he shows first that the leaves of the lamination $\sigma$ are minimal
planes. Then by using topological arguments, Gabai proves that these planes are not only minimal, but also area
minimizing. Then, he shows that this lamination must stay in a neighborhood of the convex hull of $\Gamma$, i.e.
$\sigma \subset N_C (CH(\Gamma))$ where $CH(\Gamma)$ is the convex hull of $\Gamma$ and $C$ is a constant independent
of $\Gamma$. Then,  he shows that $\PI \sigma = \Gamma$ and finishes the proof.
\end{pf}

\begin{rmk} Until this paper, all the existence results on this problem came out via the techniques of geometric measure
theory. The disadvantage of geometric measure theory is that it is very powerful with absolutely area minimizing
submanifolds, but you have to work very hard to get results in fixed topological type case. On the other hand, Gabai's
techniques are very natural for the fixed topological type case as you can control the limiting process and limiting
object topologically.
\end{rmk}

Later, the author generalized Gabai's results to the Gromov hyperbolic $3$-spaces with cocompact
metric.

\begin{thm} \cite{Co1}
Let $X$ be a Gromov hyperbolic $3$-space with cocompact metric. Let $\Gamma$ be a simple closed
curve in $S^2_\infty(X)$. Then, there exists a $D^2$-limit lamination $\sigma$ whose leaves are
area minimizing planes in $X$ with $\PI \sigma =\Gamma$.
\end{thm}

\section{Boundary Regularity at Infinity}

After the above existence theorems, the next natural question was the regularity of the hypersurfaces $\Sigma$ obtained
as a solution of the asymptotic Plateau problem. By the interior regularity theorems of geometric measure theory,
$\Sigma$ is real analytic hypersurface of $\BHH$ away from a singular subset of Hausdorff dimension $n-7$. The question
is the behavior of the hypersurfaces near infinity, i.e. the boundary regularity at infinity. In other words, if
$\Sigma$ is an area minimizing hypersurface in $\BHH$, then what can be said about the boundary regularity of
$\overline{\Sigma}$ in $\overline{\BHH}$?

The first main result about this problem came from Hardt and Lin in \cite{HL}. By using geometric measure theory
methods, they showed that near infinity, $\overline{\Sigma}$ is as regular as the asymptotic boundary for
$C^{1,\alpha}$ asymptotic boundary data.

\begin{thm} \cite{HL}
Let $\Gamma$ be a $C^{1,\alpha}$ codimension-1 submanifold of $\SI$ where $0<\alpha\leq 1$. If
 $\Sigma$ is a complete, absolutely area minimizing locally integral $n$-current in $\BHH$ with $\PI
\Sigma = \Gamma$. Then, near $\Gamma$, $\Sigma\cup\Gamma$ is union of $C^{1,\alpha}$ submanifolds
with boundary in Euclidean metric on $\overline{\BHH}$. These submanifolds have disjoint analytic
interiors, and they meet $\SI$ orthogonally at $\Gamma$.
\end{thm}

Also, if we take the upper half space model for $\BHH$, then $\BR^{n}\times \{0\}\cup\{\infty\}$ would represent the
asymptotic sphere. Then, for a given $C^1$ hypersurface $\Gamma$ in $\BR^{n}\times \{0\}$, there is $\rho_\Gamma$ with
$(\Sigma\cup\Gamma) \cap \{y<\rho_\Gamma\}$ is a finite union of $C^1$ submanifolds with boundary which can be viewed
as a graph over $\Gamma \times [0,\rho_\Gamma)$.

This result is very interesting as an area minimizing hypersurface in $\BHH$ has better regularity near asymptotic
boundary than in the interior. In other words, if $\Sigma$ is an area minimizing hypersurface in $\BHH$ with $\PI
\Sigma = \Gamma$ as above, $\Sigma$ might have a singular set of Hausdorff dimension $n-7$, but this set must stay in
the bounded part of $\Sigma$ as $(\Sigma\cup\Gamma) \cap \{y<\rho_\Gamma\}$ is a finite union of $C^1$ submanifolds
with boundary. In order to get this result, Hardt and Lin first get an interior regularity result "near infinity" by
showing that $\Sigma$ can be expressed as a union of graphs of finitely many analytic functions on vertical planes
tangent to $\Gamma$. Then by using this interior regularity "near infinity" result, and hyperbolic isometries, they
prove the regularity at boundary. In particular, if there was a sequence of singular points escaping to infinity (or
converging to a point in asymptotic boundary), by rescaling $\Sigma$ with hyperbolic isometries, they get new area
minimizing hypersurfaces, and the images of the singular points in these new area minimizing hypersurfaces would
contradict the earlier interior regularity result.

Later, by studying the following quasilinear, non-uniformly elliptic equation whose solutions are minimal hypersurfaces
in hyperbolic space, Lin and Tonegawa got higher regularity near asymptotic boundary. In the upper half space model of
$\BHH$, let $\Omega\subset \BR^n\times\{0\}$ be a domain and $f:\Omega\to \BR^+$ be a function. Consider
$graph(f)=\Sigma_f$ which defines a hypersurface in $\BHH$. The volume of $\Sigma^K_f=\Sigma_f\cap\{K\times \BR^+\}$
where $K$ is a compact subset of $\Omega$ can be described as follows:

$$vol(\Sigma^K_f) = \int_K f^{-n}\sqrt{1+|\nabla f|^2} dx$$

Then, the corresponding Euler-Lagrange equation of this variational integral would give the following Dirichlet
problem:

$$\begin{array}{ll}
\nabla f- \frac{f_i.f_j}{1+|df|^2}f_{ij}+\frac{n}{f}=0 & \mbox{in} \ \Omega\\
f>0 & \mbox{in}\  \Omega\\
f=0 & \mbox{in} \ \partial \Omega
\end{array}$$

\noindent where $|df|^2=\sum_{i=1}^n f_i^2$. In \cite{An2}, Anderson showed the existence and the uniqueness of the
solution to this Dirichlet problem provided that $\partial \Omega$ has nonnegative mean curvature with respect to
inward normal in $\BR^n\times\{0\}$.

If one wants to focus on the boundary regularity of the solution of this Dirichlet problem, an equivalent local
description of the problem can be given by considering $graph(f)$ near a point of the asymptotic boundary as a graph
over a vertical plane which is tangent to the asymptotic boundary at the given point. In other words, let
$\Gamma=\partial \Omega$ be at least $C^1$. Let $P$ be the vertical tangent plane to $\Gamma$ at $p$. By using
hyperbolic isometries, we can assume $p=0$ in $\BR^n\times \{0\}$ and $P$ is the plane $\{(x,0,y) \in \overline{\BHH} \
 | \ (x,0)\in \BR^n \ \mbox{and} \ y\geq 0\}$. Then after scaling with hyperbolic isometries, we can formulate the
problem as follows: Let $u: D \to \BR$ where $D=\{(x,0,y) \in P \ | \ |x|\leq 1 \ \mbox{and} \ 0\leq y \leq 1 \}$

$$\begin{array}{ll}
y(\nabla u- \frac{u_i.u_j}{1+|du|^2}u_{ij})-n.u_y=0 & \mbox{in} \ D\\
u(x,0,0)=\varphi(x)
\end{array}$$

\noindent where $u(x,0,0)=\varphi(x)$ is given by $\Gamma$ near $p$. Hence the question becomes whether $u$ is as
smooth as $\varphi$.

Lin studied first this quasi-linear degenerate elliptic partial differential equation in \cite{Li1} and got the
following result.

\begin{thm} \cite{Li1}
Let $\Gamma$ be a $C^{k,\alpha}$ codimension-1 submanifold of $\SI$ where $1\leq k\leq n-1$ and $0\leq\alpha\leq 1$ or
$k=n$ and $0\leq\alpha<1$. If $\Sigma$ is a complete area minimizing hypersurface in $\BHH$ with $\PI \Sigma = \Gamma$.
Then, near $\Gamma$, $\Sigma\cup\Gamma$ is union of $C^{k,\alpha}$ submanifolds with boundary in Euclidean metric on
$\overline{\BHH}$.
\end{thm}

Later, Tonegawa completed Lin's results for higher regularity case by studying further the above PDE, and finished off
the problem by giving the following very interesting parity in \cite{To}.

\begin{thm} \cite{To}
Let $\Gamma$ be a $C^{k,\alpha}$ codimension-1 submanifold of $\SI$ and $\Sigma$ be a complete area minimizing
hypersurface in $\BHH$ with $\PI \Sigma = \Gamma$. Let $k\geq n+1$ and $0<\alpha< 1$. Then,

1. If $n$ is even, then $\Sigma \cup\Gamma$ is a $C^{k,\alpha}$ submanifold with boundary near $\Gamma$.

2. If $n$ is odd, then $\Sigma\cup \Gamma$ may not be a $C^{n+1}$ submanifold with boundary near $\Gamma$ in general.
\end{thm}

This is a very interesting result as it gives a very subtle relation between the dimension and the asymptotic
regularity of area minimizing hypersurfaces. In particular, for $n$ odd, Tonegawa gives a necessary and sufficient
condition that $\Gamma$ has to satisfy in the form of a PDE in order to recover $C^{k,\alpha}$ regularity. Hence, when
$n$ is odd, if $\Gamma$ does not satisfy this PDE, $\Sigma\cup\Gamma$ cannot be smoother than $C^{n+1}$ even though
$\Gamma$ is very smooth. Note also that in \cite{To}, Tonegawa studied a more general form of the PDE above and
generalized these results to Constant Mean Curvature (CMC) hypersurfaces in $\BHH$ (See Section 6.2).

For the higher codimension case ($k<n$), by the interior regularity results of geometric measure theory, the absolutely
area minimizing k-currents are smoothly embedded $k$-submanifolds in $\BHH$ except for a singular set of Hausdorff
dimension $k-2$. For the boundary regularity at infinity in this case, Lin also showed the {\em existence} of an area
minimizing $k$-current $\Sigma$ in $\BHH$ which is as regular as the boundary at infinity, where $\Gamma= \PI \Sigma$
is a $C^{1,\alpha}$ smooth closed $k-1$-submanifold in $\SI$.

\begin{thm} \cite{Li2}
Let $\Gamma$ be a $C^{1,\alpha}$ smooth closed $k-1$-submanifold in $\SI$. Then there exists a complete area minimizing
$k$-current in $\BHH$ with $\PI \Sigma =\Gamma$ such that near $\Gamma$, $\Sigma\cup\Gamma$ is a $C^{1,\alpha}$
submanifold with boundary in Euclidean metric on $\overline{\BHH}$.
\end{thm}

Note that unlike the codimension-1 case, this higher codimension case does not say {\em any} area minimizing
$k$-current with asymptotic smooth asymptotic boundary is boundary regular at infinity. This result only says the
existence of such an area minimizing current for any given smooth asymptotic data.

\section{Number of Solutions}

There are basically $3$ types of results on the number of solutions to the asymptotic Plateau problem. The first type
is the uniqueness results which classifies the asymptotic data with the unique solution to the asymptotic Plateau
problem. The next type is the generic uniqueness and generic finiteness results which came out recently. The last type
can be called as the nonuniqueness results which constructs the asymptotic data with more than one solution to the
problem.

\subsection{Uniqueness and Finiteness Results}\ \\

Next to the existence theorems, Anderson gave very interesting uniqueness and nonuniqueness results
on minimal surfaces in $\BH^3$ and area minimizing hypersurfaces in $\BHH$ in \cite{An1} and
\cite{An2}. Before visiting nonuniqueness results, we will list the uniqueness results about the
the asymptotic Plateau problem.

First, in \cite{An1}, Anderson showed that if the given asymptotic boundary $\Gamma_0$ which is a
hypersurface bounding a convex domain in $\SI$, then there exists a unique absolutely area
minimizing hypersurface $\Sigma_0$ in $\BHH$.

\begin{thm} \cite{An1}
Let $\Gamma_0$ be a hypersurface bounding a convex domain in $\SI$. Then, there exists a unique
absolutely area minimizing hypersurface $\Sigma_0$ in $\BHH$ with $\PI \Sigma_0 =\Gamma_0$.
\end{thm}

\begin{pf} (Sketch) Let $\Gamma_0$ be codimension-1 submanifold bounding a convex domain in $\SI$, and $\Sigma_0$ be an area
minimizing hypersurface in $\BHH$ with $\PI \Sigma_0 = \Gamma_0$ (Existence of $\Sigma_0$ is guaranteed by Theorem
3.1). As $\Gamma_0$ bounds a convex domain in $\SI$, we can find a continuous family of isometries $\{\varphi_t\}$ of
$\BHH$ such that $\varphi_t(\Gamma_0)=\Gamma_t$ where $\{\Gamma_t\}$ foliates whole $\SI$. Similarly, if
$\varphi_t(\Sigma_0)=\Sigma_t$, then $\PI \Sigma_t = \Gamma_t$, and as $\{\Sigma_t\}$ images of continuous family of
isometries, it foliates whole $\BHH$.

Hence, if there are two minimal hypersurfaces $M_1, M_2$ with $\PI M_i=\Gamma_0$, one of them (say $M_2$) is not a leaf
of the foliation, and $M_2$ must intersect a leaf $\Sigma_{t_0}$ of the foliation tangentially and by lying in one
side.This contradicts to the maximum principle for minimal hypersurfaces.
\end{pf}

Later, by using similar ideas, Hardt and Lin generalized this result to the codimension-1
submanifolds bounding star shaped domains in $\SI$ in \cite{HL}.

\begin{thm} \cite{HL}
Let $\Gamma_0$ be a hypersurface bounding a star shaped domain in $\SI$. Then, there exists a
unique absolutely area minimizing hypersurface $\Sigma_0$ in $\BHH$ with $\PI \Sigma_0 =\Gamma_0$.
\end{thm}

While these are the only known results on the number of solutions of the asymptotic Plateau problem
for a long time, many generic uniqueness results have come out recently in both general case and
fixed topological type case.

For the general case, the author showed that the space of closed codimension-1 submanifolds
$\Gamma$ in $\SI$ bounding a unique absolutely area minimizing hypersurface $\Sigma$ in $\BH^n$ is
dense in the space of all closed codimension-1 submanifolds in $\SI$ by using a simple topological
argument.

\begin{thm} \cite{Co7}
Let $B$ be the space of connected closed codimension-$1$ submanifolds of $\SI$, and let $B'\subset
B$ be the subspace containing the closed submanifolds of $\SI$ bounding a unique absolutely area
minimizing hypersurface in $\BH^n$. Then $B'$ is dense in $B$.
\end{thm}

\begin{pf} (Sketch) For simplicity, we will focus on the area minimizing planes in $\BH^3$. The general case is similar. Let $\Gamma_0$ be a simple
closed curve in $\Si$. First, by using Meeks-Yau exchange roundoff trick, the author establishes that if $\Gamma_1$ and
$\Gamma_2$ are two disjoint simple closed curves in $\Si$, and $\Sigma_1$ and $\Sigma_2$ are area minimizing planes in
$\BH^3$ with $\PI \Sigma_i = \Gamma_i$, then $\Sigma_1$ and $\Sigma_2$ are disjoint, too. Then, by using this result,
he shows that for any simle closed curve $\Gamma$ in $\Si$ either there exists a unique area minimizing plane $\Sigma$
in $\BH^3$ with $\PI\Sigma=\Gamma$, or there exist two {\em disjoint} area minimizing planes $\Sigma^+ , \Sigma^-$ in
$\BH^3$ with $\PI\Sigma^\pm=\Gamma_0$.

Then, take a small neighborhood $N(\Gamma_0)\subset \Si$ which is an annulus, and foliate $N(\Gamma_0)$ by simple
closed curves $\{\Gamma_t\}$ where $t\in(-\epsilon, \epsilon)$, i.e. $N(\Gamma_0) \simeq \Gamma\times (-\epsilon,
\epsilon)$. By the above fact, for any $\Gamma_t$ either there exists a unique area minimizing plane $\Sigma_t$, or
there are two area minimizing planes $\Sigma_t^\pm$ disjoint from each other. As disjoint asymptotic boundary implies
disjoint area minimizing planes, if $t_1<t_2$, then $\Sigma_{t_1}$ is disjoint and \textit{below} $\Sigma_{t_2}$ in
$\BH^3$. Consider this collection of area minimizing planes. Note that for curves $\Gamma_t$ bounding more than one
area minimizing plane, we have a canonical region $N_t$ in $\BH^3$ between the disjoint area minimizing planes
$\Sigma_t^\pm$.

\begin{figure}[t]
\mbox{\vbox{\epsfxsize=3in \epsfbox{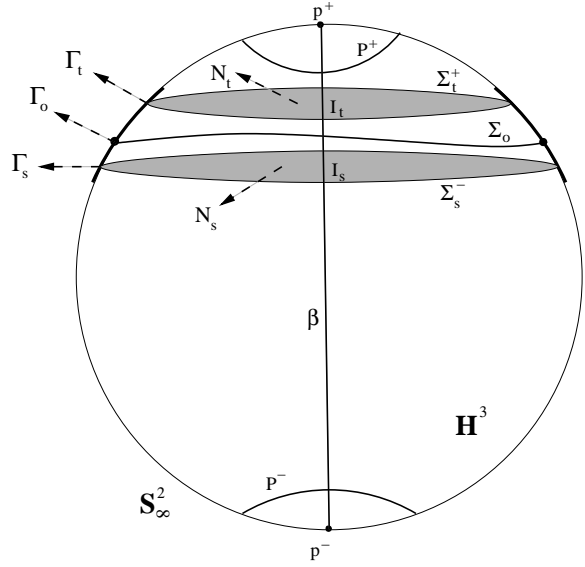}}} \caption{\label{fig:figure1} {A finite segment of geodesic $\gamma$
intersects the collection of area minimizing planes $\Sigma_t$ in $\BH^n$ asymptotic to $\Gamma_t$ in $\Si$.}}
\end{figure}

Now, the idea is to consider the {\em thickness} of the neighborhoods $N_t$ assigned to the asymptotic curves
$\{\Gamma_t\}$. Let $s_t$ be the length of the segment $I_t$ of $\beta$ (a fixed finite length transversal curve to the
collection) between $\Sigma_t^+$ and $\Sigma_t^-$, which is the {\em width} of $N_t$ assigned to $\Gamma_t$. Then, the
curves $\Gamma_t$ bounding more than one area minimizing planes have positive width, and contributes to the total
thickness of the collection, and the curves bounding unique area minimizing plane has $0$ width and do not contribute
to the total thickness. Since $\sum_{t\in(-\epsilon, \epsilon)} s_t < C$, the total thickness is finite. This implies
for only countably many $t\in(-\epsilon, \epsilon)$, $s_t>0$, i.e. $\Gamma_t$ bounds more than one area minimizing
plane. For the remaining uncountably many $t\in(-\epsilon, \epsilon)$, $s_t=0$, and there exists a unique area
minimizing plane for those $t$. This proves the space of Jordan curves of uniqueness is dense in the space of Jordan
curves in $\Si$. Then, the author shows that this space is not only dense, but also generic. Also, this technique is
quite general, and it can be generalized to many different settings \cite{Co4}, \cite{Co9}.
\end{pf}

On the other hand, there has been important progress on the number of solutions to the asymptotic Plateau problem in
fixed topological type case. Recently in \cite{Co2}, the author showed a generic finiteness result for $C^3$ smooth
Jordan curves in $\Si$ for area minimizing planes in $\BH^3$ by using geometric analysis and global analysis methods.
Later in \cite{Co3}, he improved this result to a generic uniqueness result.

\begin{thm} \cite{Co3}
Let $A$ be the space of $C^3$ simple closed curves in $\Si$. Then there exists an open dense subset
$A' \subset A$ such that for any $\Gamma \in A'$, there exists a unique area minimizing plane
$\Sigma$ with $\PI \Sigma=\Gamma$.
\end{thm}

\begin{pf} (Sketch) In \cite{Co2}, by generalizing Tomi and Tromba's global analytic techniques in \cite{TT} to hyperbolic setting,
and by using Li and Tam's powerful results \cite{LT1} and \cite{LT2}, the author showed that the boundary restriction
map $\pi$ from the space of minimal maps from $D^2$ to $\BH^3$ with $C^3$ asymptotic data to the space of the $C^3$
immersions of $S^1$ into $\Si$ is Fredholm of index $0$. Hence, the derivative of $\pi$ is isomorphism for the generic
curves.

Fix a generic curve $\Gamma$ in $\Si$. By using the inverse function theorem, there is a neighborhood $U_\Sigma$ of a
area minimizing plane $\Sigma$ in $\pi^{-1}(\Gamma)$, mapping homeomorphically into a neighborhood $V_\Gamma$ of
$\Gamma$. By taking a path $\alpha$ in $V_\Gamma$, and by considering the corresponding path $\pi^{-1}(\alpha)$ in
$U_\Sigma$, one can get a continuous family of minimal planes with disjoint asymptotic boundaries around $\Sigma$.
Then, the author shows that this continuous family of minimal planes is indeed a foliation by area minimizing planes of
a neighborhood of $\Sigma$. This implies the uniqueness of the area minimizing plane in $\BH^3$ spanning $\Gamma$. Then
the author proves that the same is true for any curve in a neighborhood of a generic curve, and gets an open dense
subset of the $C^3$ Jordan curves in $\Si$ with the uniqueness result.
\end{pf}

Recently, in \cite{AM}, Alexakis and Mazzeo generalized this result to any surface of genus $k$ by using different
methods. In \cite{Co2}, the author works with the space of parametrizations of minimal planes in $\BH^3$, and hence, in
order to get a generic finiteness result, he needs to deal with different parametrizations of the same minimal plane.
In \cite{AM}, Alexakis and Mazzeo showed that if $\mathcal{M}_k$ is the moduli space of all complete minimal surfaces
of genus $k$ in $\BH^3$ with asymptotic boundary curve a $C^{3,\alpha}$ simple closed curve in $\Si$, and $\xi$ is the
space of $C^{3,\alpha}$ simple closed curve in $\Si$, then the boundary restriction map $\pi_k:\mathcal{M}_k\to \xi$ is
Fredholm of index $0$ (see also Section 7.2). Moreover, they also showed that $\pi_k$ is not only Fredholm of index
$0$, but also proper (Theorem 4.3 in \cite{AM}). Hence by Sard-Smale theorem \cite{Sm}, this implies a generic
finiteness result for minimal surfaces of genus $k$. In other words, for a generic $C^{3,\alpha}$ simple closed curve
$\Gamma$ in $\Si$, there exist finitely many complete minimal surfaces $\Sigma$ of genus $k$ in $\BH^3$ with $\PI
\Sigma = \Gamma$. Indeed, their result also applies to convex cocompact hyperbolic $3$-manifolds, too.

Note that the above generic uniqueness result for area minimizing planes requires some smoothness condition on the
curves. Later, the author improved his result by removing the smoothness condition. This time, the author uses
topological methods instead of techniques of global analysis. The technique is essentially same with the area
minimizing hypersurfaces case mentioned above.

\begin{thm} \cite{Co7}
Let $A$ be the space of simple closed curves in $\Si$ and let $A'\subset A$ be the subspace
containing the simple closed curves in $\Si$ bounding a unique area minimizing plane in $\BH^3$.
Then, $A'$ is generic in $A$, i.e. $A-A'$ is a set of first category.
\end{thm}

\begin{rmk} Note that the same result is true for area minimizing surfaces in $\BH^3$, too \cite{Co7}.
By these results, the asymptotic Plateau problem generically has a unique solution in both area minimizing surfaces in $\BH^3$ case and area
minimizing planes in $\BH^3$ case. In higher dimensions, the closed codimension-1 submanifolds in $\SI$ bounding a
unique absolutely area minimizing hypersurface in $\BHH$ are only dense in the closed codimension-1 submanifolds in
$\SI$. However, by using the similar ideas, by fixing the topological type of the closed codimension-1 submanifold in
$\SI$, it might be possible to get some generic uniqueness result, too.
\end{rmk}

\begin{rmk} Notice that except the convex and star-shaped asymptotic boundary cases, all the uniqueness results on the asymptotic Plateau
are about area minimizing surfaces or area minimizing planes. Unfortunately, the techniques used for these results
cannot be extended to the minimal surfaces or minimal planes cases. The main obstacle here is that while two area
minimizing surfaces with disjoint asymptotic boundaries must be disjoint, the same statement may not be true for
minimal surfaces. In any case, it would be interesting problem to study this case in order to understand whether the
simple closed curves in $\Si$ bounding a unique minimal surface (or plane) is dense in the space of simple closed
curves in $\Si$ or not. The author believes that the similar statements are not true in minimal surfaces (or planes)
case.
\end{rmk}

\subsection{Nonuniqueness Results}\ \\

Beside his existence results, Anderson also gave many different nonuniqueness results for the
asymptotic Plateau problem in the fixed topological type in \cite{An2}.

\begin{thm} \cite{An2}
There exists a simple closed curve $\Gamma$ in $\Si$ such that there are infinitely many complete
minimal surfaces $\{\Sigma_i\}$ in $\BH^3$ with $\PI \Sigma_i = \Gamma$.
\end{thm}

For the proof of this theorem, Anderson first constructs a simple closed curve such that the
absolutely area minimizing surface given by his existence theorem is not a plane (positive genus)
(a similar construction can be found in \cite{Ha}). Then, by modifying this curve, he constructs a
curve $\Gamma$ with the same property such that it is also a limit set for a quasi-Fuchsian group
$\Lambda$. Since the absolutely area minimizing surface $\Sigma$ is $\Lambda$ invariant and has
positive genus, this implies the area minimizing surface $\Sigma / \Lambda$ in the compact
hyperbolic manifold $\BH^3$ is not $\pi_1$-injective. This implies that the absolutely area
minimizing surface $\Sigma$ in $\BH^3$ with $\PI \Sigma = \Gamma$ must have infinite genus. Then,
by using this property, he shows that there exist infinitely many complete {\em minimal} surfaces
asymptotic to $\Gamma$.

Note that this result shows nonuniqueness for \textit{minimal surfaces} for fixed topological type. Later, the author
show nonuniqueness for area minimizing planes (surfaces) in $\BH^3$ case. In particular, the author shows that there
are simple closed curves in $\Si$ bounding more than one area minimizing plane (surface).

\begin{thm} \cite{Co7}
There exists a simple closed curve $\Gamma$ in $\Si$ such that there are more than one area minimizing plane (surface)
$\{\Sigma_i\}$ in $\BH^3$ with $\PI \Sigma_i = \Gamma$.
\end{thm}

\begin{rmk} In the nonuniqueness results above, only Hass' result gives an explicit example of a simple closed curve in
$\Si$ bounding more than one minimal surfaces in $\BH^3$. All other results on nonuniqueness so far shows the existence
of such a curve, but it does not give one. So, it would be interesting to construct an explicit simple closed curve in
$\Si$ bounding more than one area minimizing surface (or plane).
\end{rmk}

\begin{rmk} Although there are many examples of simple closed curves in $\Si$ bounding more than one minimal surface or more than
one area minimizing surface (or plane) in $\BH^3$, there is no example in higher dimensions so far. It would be
interesting to generalize the nonuniqueness results to higher dimensions by showing whether there exists a closed
codimension-1 submanifold in $\SI$ bounding more than one absolutely area minimizing hypersurfaces in $\BHH$.
\end{rmk}

\section{CMC Hypersurfaces}

After many important results on minimal hypersurfaces in hyperbolic space, like existence,
regularity, etc., the question of generalization of these results to constant mean curvature (CMC)
hypersurfaces was naturally raised: For a given codimension-1 submanifold $\Gamma$ in $\SI$, does
there exists a complete CMC hypersurface $\Sigma$ with specified mean curvature $H$ in $\BHH$ and
$\PI \Sigma = \Gamma$?

For simplicity, from now on, we will call CMC hypersurfaces with mean curvature $H$ as \textit{$H$-hypersurfaces}.

Note that for this generalization of the asymptotic Plateau problem, we need to assume that $|H|<1$ (after fixing an
orientation on $\BHH$). This is because it is impossible to have a complete  $H$-hypersurface $\Sigma$ in $\BHH$ with
$|H|\geq 1$ and $\PI \Sigma = \Gamma$ as we can always find a horosphere ($H=1$) in $\BHH$ with tangential intersection
with such a $\Sigma$ which contradicts to the maximum principle.

We should also note that $H$-hypersurfaces in $\BH^3$ with $H=1$ and $H>1$ are also an area of active research. A basic
reference for CMC hypersurfaces in hyperbolic space with $H>1$ would be \cite{KKMS}. For the case $H=1$, we refer to
Rosenberg's survey \cite{Ro}, and Bryant's seminal paper \cite{Br} where he showed that any minimal surface in $\BR^3$
is isometric to a CMC surface in $\BH^3$ with $H=1$.

We should point out that the generalization of area minimizing (or minimal) hypersurfaces to CMC
hypersurfaces is quite natural. As we see the minimal hypersurfaces ($H=0$) as the critical points
of the area functional, CMC hypersurfaces occurs as the critical points of some modification of the
area functional with a volume constraint. In particular, let $\Sigma^n$ be a compact hypersurface,
bounding a domain $\Omega^{n+1}$ in some ambient Riemannian manifold. Let $A$ be the area of
$\Sigma$, and $V$ be the volume of $\Omega$. Let's vary $\Sigma$ through a one parameter family
$\Sigma_t$, with corresponding area $A(t)$ and volume $V(t)$. If $f$ is the normal component of the
variation, and $H$ is the mean curvature of $\Sigma$, then we get $A'(0) = -\int_\Sigma n H f$, and
$V'(0)=\int_\Sigma f$ where $n$ is the dimension of $\Sigma$, and $H$ is the mean curvature.

Now, let $\Sigma$ be a hypersurface with boundary $\Gamma$. We fix a hypersurface $M$ with $\partial M = \Gamma$, and
define $V(t)$ to be the volume of the domain bounded by $M$ and $\Sigma_t$. Now, we define a new functional as a
combination of $A$ and $V$. Let $I_H(t)= A(t) + n H V(t)$. Note that $I_0(t)=A(t)$. If $\Sigma$ is a critical point of
the functional $I_H$ for any variation $f$, then this will imply $\Sigma$ has constant mean curvature $H$ \cite{Gu}.
Note that critical point of the functional $I_H$ is independent of the choice of the hypersurface $M$ since if
$\widehat{I}_H$ is the functional which is defined with a different hypersurface $\widehat{M}$, then $I_H -
\widehat{I}_H = C$ for some constant $C$. On the other hand, we will call $\Sigma$ a \textit{minimizing CMC
hypersurface} if $\Sigma$ is the absolute minimum of the functional $I_H$ among hypersurfaces with the same boundary.
From this point of view, CMC hypersurfaces are natural generalization of minimal hypersurfaces and area minimizing
hypersurfaces as the area functional is just $H=0$ case for the functional $I_H$. This point of view is very useful and
essential to generalize the geometric measure theory methods developed for area minimizing case to CMC case as in
\cite{To} and \cite{AR}.

Now, we continue with the basic notions on $H$-hypersurfaces in $\BHH$. We fix a codimension-$1$ closed submanifold
$\Gamma$ in $\SI$. $\Gamma$ separates $\SI$ into two parts, say $\Omega^+$ and $\Omega^-$. By using these domains, we
will give orientation to hypersurfaces in $\BHH$ asymptotic to $\Gamma$. With this orientation, mean curvature $H$ is
positive if the mean curvature vector points towards positive side of the hypersurface, negative otherwise. The
following fact is known as maximum principle.

\begin{lem} $[$Maximum Principle$]$
Let $\Sigma_1$ and $\Sigma_2$ be two hypersurfaces in a Riemannian manifold which intersect at a common point
tangentially. If $\Sigma_2$ lies in positive side (mean curvature vector direction) of $\Sigma_1$ around the common
point, then $H_1$ is less than or equal to $H_2$ ($H_1 \leq H_2$) where $H_i$ is the mean curvature of $\Sigma_i$ at
the common point. If they do not coincide in a neighborhood of the common point, then $H_1$ is strictly less than $H_2$
($H_1<H_2$).
\end{lem}

The other important notion about CMC Hypersurfaces in $\BHH$ is the generalization of the convex hull property to this
context. Now, let $\Gamma$ be a codimension-$1$ submanifold of $\SI$ and orient all spheres accordingly. If $T$ is a
round $n-1$-sphere in $\SI$, then there is a unique $H$-hypersurface $P_H$ in $\BHH$ asymptotic to $T$ for $-1<H<1$
\cite{NS}. $T$ separates $\SI$ into two parts $\Delta^+$ and $\Delta^-$. Similarly, $P_H$ divides $\BHH$ into two
domains $D_H^+$ and $D_H^-$ with $\PI D_H^\pm = \Delta^\pm$. We will call these regions as \textit{$H$-shifted
halfspaces}. If the asymptotic boundary of a $H$-shifted halfspace contains $\Gamma$, then we will call this
$H$-shifted halfspace as \textit{supporting $H$-shifted halfspace}. i.e. if $A\subset\Delta^+$, then $D_H^+$ is a
supporting $H$-shifted halfspace. Then the {\em $H$-shifted convex hull} of $\Gamma$, $CH_H(\Gamma)$ is defined as the
intersection of all supporting closed $H$-shifted halfspaces of $\BHH$.

Now, the generalization of convex hull property of minimal hypersurfaces in $\BHH$ to $H$-hypersurfaces in $\BHH$ is as
follows \cite{Co4}. Similar versions of this result have been proved by Alencar-Rosenberg in \cite{AR}, and by Tonegawa
in \cite{To}.

\begin{lem} \cite{To}, \cite{AR}, \cite{Co1}
Let $\Sigma$ be a $H$-hypersurface in $\BHH$ where $\PI\Sigma = \Gamma$ and $|H|<1$. Then $\Sigma$ is in the
$H$-shifted convex hull of $\Gamma$, i.e. $\Sigma \subset CH_H(\Gamma)$.
\end{lem}

\subsection{Existence}\ \\

In the following decade after Anderson's existence (\cite{An1}, \cite{An2}) and Hardt-Lin's
regularity results (\cite{HL},\cite{Li1}), there have been many important generalizations of these
results to CMC hypersurfaces in hyperbolic space. In \cite{To}, Tonegawa generalized Anderson's
techniques to this case, and proved existence for CMC hypersurfaces by using geometric measure
theory methods. In the same year, by using similar techniques, Alencar and Rosenberg got a similar
existence result in \cite{AR}.

\begin{thm} \cite{To}, \cite{AR}
Let $\Gamma$ be a codimension-$1$ closed submanifold in $\SI$, and let $|H|<1$. Then there exists a
CMC hypersurface $\Sigma$ with mean curvature $H$ in $\BHH$ where $\PI \Sigma = \Gamma$. Moreover,
any such CMC hypersurface is smooth except a closed singularity set of dimension at most $n-7$.
\end{thm}

We should also note that Nelli and Spruck showed existence of a CMC hypersurface asymptotic to
$C^{2,\alpha}$ codimension-1 submanifold $\Gamma$ which is the boundary of a mean convex domain in
$\SI$ by using analytic techniques in \cite{NS}. Later, Guan and Spruck generalized this result to
$C^{1,1}$ codimension-1 submanifolds bounding star shaped domains in $\SI$.

\begin{thm} \cite{GS}
Let $\Omega$ be a star shaped (mean convex in \cite{NS}) domain in $\SI$ where $\Gamma=\partial
\Omega$ is $C^{1,1}$ ($C^{2,\alpha}$ in \cite{NS}) codimension-1 submanifold in $\SI$. Then, for
any $0<H<1$, there exists a complete smoothly embedded CMC hypersurface $\Sigma$ with mean
curvature $H$ and $\PI \Sigma = \Gamma$. Moreover, $\Sigma$ can be represented as a graph of a
function $u \in C^{1,1}(\overline{\Omega})$ ($u \in C^{2,\alpha}(\overline{\Omega})$ in \cite{NS}).
\end{thm}

Even though this second existence result is for fairly restricted class of asymptotic boundary data
(star shaped condition), the CMC hypersurfaces obtained are smoothly embedded with no singularity
in any dimension (unlike the first one), and they can be represented as a graph like $x_{n+1} = u$
for a function $u \in C^{1,1}(\overline{\Omega})$ in half space model for $\BHH$. We should also
note that, in \cite{AA}, Aiyama and Akutagawa gave a completely different construction for CMC
surfaces of disk type in $\BH^3$ with asymptotic boundary $C^{1,\alpha}$ smooth simple closed curve
in $\Si$ by studying a Dirichlet problem at infinity.

\subsection{Boundary Regularity at Infinity}\ \\

Beside the existence results, in \cite{To}, Tonegawa studied the following quasi-linear degenerate elliptic PDE which
is a more general form of the PDE in Section 4 for $H$-hypersurfaces with $|H|<1$, and got important regularity results
for these hypersurfaces near asymptotic boundary.

$$\begin{array}{ll}
y(\nabla u- \frac{u_i.u_j}{1+|du|^2}u_{ij})-n(u_y-H\sqrt{1+|Du|^2})=0 & \mbox{in} \ D\\
u(x,0,0)=\varphi(x)
\end{array}$$

For $k\leq n$, Tonegawa generalized the Lin's result for minimal hypersurfaces ($H=0)$ in \cite{Li1}.

\begin{thm} \cite{To}
Let $\Gamma$ be $C^{k,\alpha}$ codimension-1 submanifold in $\SI$ where $1\leq k  \leq n-1$ and
$0\leq \alpha \leq 1$ or $k=n$ and $0\leq \alpha <1$. If $\Sigma$ is a complete CMC hypersurface in
$\BHH$ with $\PI \Sigma = \Gamma$, then $\Sigma\cup\Gamma$ is a $C^{k,\alpha}$ submanifold with
boundary in $\overline{\BHH}$ near $\Gamma$.
\end{thm}

On the other hand, Tonegawa showed that for higher regularity case, $H=0$ case is fairly different form the $H\neq 0$
case. As we mentioned in Section 4, in $H=0$ case, Tonegawa showed that when $n$ is even the higher regularity is
always true, but when $n$ is odd, the higher regularity depends on the asymptotic boundary $\Gamma$ (Theorem 4.3). In
the $H\neq 0$ case, Tonegawa got a very surprising result that while the similar result is true for $n=2$, it is not
true for $n=4$.

\begin{thm} \cite{To}
{\bf a. ($n=2$ case)} Let $\Gamma$ be a $C^{k,\alpha}$ smooth simple closed curve in $\Si$ with $k\geq n+1=3$,
$0<\alpha <1$. Let $\Sigma$ be a $H$-hypersurface in $\BH^3$ with $\PI \Sigma = \Gamma$. Then, $\Sigma \cup \Gamma$ is
a $C^{k,\alpha}$ submanifold with boundary near $\Gamma$.\\

\noindent {\bf b. ($n=4$ case)} For $n=4$, $H\neq 0$ and $|H|<1$, there exists a smooth $\Gamma$ such that $\Sigma \cup
\Gamma$ is not a $C^{n+1}=C^5$ submanifold with boundary where $\Sigma$ is a $H$-hypersurface with $\PI \Sigma =
\Gamma$.
\end{thm}

We should also note that by studying the PDE above, or by using some barrier arguments, it is not hard to show that the
intersection angle $\theta_H$ between an $H$-hypersurface and the asymptotic boundary $\SI$ is
$\arctan(\frac{\sqrt{1-H^2}}{H})$ \cite{To}. In other words, let $\Gamma$ be a codimension-1 submanifold in $\SI$, and
$\Sigma$ be a $H$-hypersurface in $\BHH$ with $\PI \Sigma = \Gamma$. Then for any $p\in \Gamma$, the angle $\theta_H$
between $\Sigma \cup \Gamma$ and $\SI$ at $p$ would be $\arctan(\frac{\sqrt{1-H^2}}{H})$.\\

\subsection{Number of Solutions}\ \\

By using analytic techniques, Nelli and Spruck generalized Anderson's uniqueness result for mean
convex domains in area minimizing hypersurfaces case to CMC context in \cite{NS}. Then, Guan and
Spruck extended Hardt and Lin's uniqueness results for star-shaped domains in area minimizing
hypersurfaces case to CMC hypersurfaces in hyperbolic space in \cite{GS}.

\begin{thm} \cite{GS}
Let $\Omega$ be a star shaped (mean convex in \cite{NS}) domain in $\SI$ where $\Gamma=\partial
\Omega$ is $C^{1,1}$ ($C^{2,\alpha}$ in \cite{NS}) codimension-1 submanifold in $\SI$. Then, for
any $0\leq H<1$, there exists a unique complete CMC hypersurface $\Sigma$ with mean curvature $H$
and $\PI \Sigma = \Gamma$.
\end{thm}

On the other hand, the author got a generic uniqueness result for CMC hypersurfaces by generalizing his methods in
\cite{Co7}. In particular, he defined the notion of {\em minimizing CMC hypersurfaces} as generalizations of area
minimizing hypersurfaces. In other words, as minimal hypersurfaces are critical points of the area functional, and area
minimizing hypersurfaces are not just critical but minimum points of the functional, the same generalization is defined
for CMC hypersurfaces in \cite{Co4}. The CMC hypersurfaces are the hypersurfaces with constant mean curvature and
corresponds to critical points of the functional $I_H(t)= A(t) + n H V(t)$, and {\em minimizing CMC hypersurfaces}
corresponds to minimizers of the functional $I_H$. Note that the existence result Theorem 6.3 by Tonegawa and
Alencar-Rosenberg indeed gives minimizing CMC hypersurfaces.

\begin{thm} \cite{Co4}
Let $A$ be the space of codimension-$1$ closed submanifolds of $\SI$, and let $A'\subset A$ be the subspace containing
the closed submanifolds of $\SI$ bounding a unique minimizing CMC hypersurface with mean curvature $H$ in $\BHH$. Then
$A'$ is generic in $A$, i.e. $A-A'$ is a set of first category.
\end{thm}

On the other hand, there is no result for nonuniqueness of CMC hypersurfaces. In particular, there is no known example
of a codimension-1 submanifold $\Gamma$ in $\SI$ such that $\Gamma$ is the asymptotic boundary of more than one CMC
hypersurface with mean curvature $H$ for any $0<H<1$. For $H=0$, Anderson \cite{An2}, Hass \cite{Ha}, and the author
\cite{Co7} gave such examples. It might be possible to generalize these techniques to prove nonuniqueness in CMC case
for any $H\in (-1,1)$.

\newpage

\subsection{Foliations of Hyperbolic Space}\ \\

While discussing the uniqueness of CMC hypersurfaces for a given asymptotic data in asymptotic boundary, there is a
related problem in the subject: For a given codimension-$1$ closed submanifold $\Gamma$ in $\SI$, does the family of
CMC hypersurfaces $\{\Sigma_H\}$ with mean curvature $H$ foliates $\BHH$ or not, where $-1<H<1$ and $\PI \Sigma_H =
\Gamma$. This question is related with uniqueness question as existence of such a foliation automatically implies the
uniqueness of CMC hypersurface $\Sigma_H$ with mean curvature $H$ where $\PI \Sigma_H = \Gamma$ by maximum principle.
In the reverse direction, the author showed the following result.

\begin{thm} \cite{Co10}
Let $\Gamma$ be a $C^{2,\alpha}$ closed codimension-1 submanifold in $\SI$. Also assume that for any $H\in (-1,1)$,
there exists a unique CMC hypersurface $\Sigma_H$ with $\PI \Sigma_H =\Gamma$. Then, the collection of CMC
hypersurfaces $\{\Sigma_H\}$ with $H\in (-1,1)$ foliates $\BHH$.
\end{thm}

\begin{pf} (Sketch) First, by using the boundary regularity results in \cite{To} and some cut-paste arguments similar
to exchange roundoff trick, the author shows that two different minimizing $H$-hypersurfaces with same asymptotic
boundary must be disjoint (See Figure 2). In particular, he proves that if $\Gamma$ is a $C^{2,\alpha}$ closed
codimension-1 submanifold in $\SI$, and $\Sigma_{H_1}$ and $\Sigma_{H_2}$ are minimizing CMC hypersurfaces in $\BHH$
with $\PI \Sigma_{H_i}=\Gamma_i$ and $-1<H_1<H_2<1$, then $\Sigma_{H_1}$ and $\Sigma_{H_2}$ are disjoint. Hence,
$\{\Sigma_H\}$ for $-1<H<1$ is a disjoint family of hypersurfaces in $\BHH$. Now, there are two points to check to show
that $\{\Sigma_H\}$ foliates $\BHH$. First point is that there is no gap between the leaves of $\{\Sigma_H\}$, and the
second point is that $\{\Sigma_H\}$ fills $\BHH$.

\begin{figure}[b]

\relabelbox  {\epsfxsize=2.7in

\centerline{\epsfbox{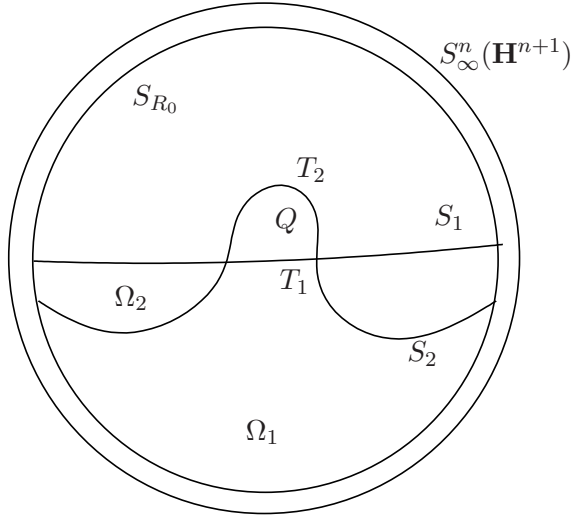}}}

\relabel{1}{$\SI$}

\relabel{2}{$S_{R_0}$}

\relabel{3}{$T_2$}

\relabel{4}{$T_1$}

\relabel{5}{$Q$}

\relabel{6}{$S_1$}

\relabel{7}{$S_2$}

\relabel{8}{$\Omega_1$}

\relabel{9}{$\Omega_2$}

\endrelabelbox

\caption{\label{fig:figure2} \small For $0<H_1<H_2<1$, $S_1$ is above $S_2$ near the boundary of the ball $B_{R_0}(p)$
by \cite{To}}.

\end{figure}

For the first point the idea is to use the assumption that $\Gamma$ bounds a unique $H$-hypersurface for any $H\in
(-1,1)$. If there was a gap between the family $\{\Sigma_H\}_{H\in(-1,H_0]}$ and $\{\Sigma_H\}_{H\in(H_0,1)}$, then
constructing a sequence of hypersurfaces $\{S_i\}$ such that $S_i\subset \Sigma_{H_i}$ where $H_i\searrow H_0$ and
$\partial S_i \to \Gamma$. Then, by showing that $S_i\to \Sigma_{H_0}'$ where $\Sigma_{H_0}'$ is another minimizing
$H_0$-hypersurface with $\PI \Sigma_{H_0}' = \Gamma$, he gets a contradiction as $\Gamma$ must bound a unique
$H_0$-hypersurface in $\BHH$.

For the second point, if $\{\Sigma_H\}$ family of hypersurfaces does not fill $\BHH$, then by constructing a suitable
horosphere in the unfilled region, and by using the maximum principle, the author gets a contradiction.
\end{pf}

Hence, by the uniqueness results in \cite{GS} and \cite{NS}, for the star shaped asymptotic data and mean convex
asymptotic data, the above result gives positive answer for the question. Note that in \cite{CV}, Chopp and Velling
studied this problem by using computational methods, and had an interesting result that for many different type of
curves in $\Si$, CMC surfaces gives a foliation of $\BH^3$.

On the other hand, recently in \cite{Wa}, Wang showed that if a quasi-Fuchsian $3$-manifold $M$ contains a minimal
surface whose principle curvature is less than $1$, than $M$ admits a foliation by CMC surfaces by using volume
preserving mean curvature flow. If we lift this foliation to the universal cover, we get a foliation of $\BH^3$ by CMC
surfaces with same asymptotic boundary $\Gamma$ where $\Gamma$ is a simple closed curve in $\Si$ and it is the limit
set of the quasi-Fuchsian $3$-manifold $M$. However, the limit set of quasi-Fuchsian manifolds are far from being
smooth, even they contain no rectifiable arcs (\cite{Be}). Existence of one smooth point in the limit set implies the
group being Fuchsian, which means the limit set is a round circle in $\Si$. Hence, in addition to smooth examples in
\cite{Co10}, \cite{Wa} gives completely nonrectifiable simple closed curve examples where CMC hypersurfaces with the
given asymptotic data foliate the hyperbolic space. Also in \cite{Wa}, Wang constructs a simple closed curve $\Gamma$
in $\Si$ (as limit set of a quasi-Fuchsian 3-manifold) which is similar to the one in \cite{Ha}, where there cannot be
a foliation of $\BH^3$ by CMC surfaces with asymptotic boundary $\Gamma$.

\section{Further Results}

Other than existence, regularity and number of solutions to the asymptotic Plateau problem, there
have been other important features which are studied.

\subsection{Properly Embeddedness}\ \\

The properly embeddedness of the solution of the asymptotic Plateau problem is one of the interesting problems which is
under investigation. Namely, the question is whether a solution to the asymptotic Plateau problem $\Sigma$ with $\PI
\Sigma = \Gamma$ where $\Gamma$ is a codimension-1 closed submanifold in $\SI$ is properly embedded, or not? In other
words, if $\varphi: S\to \BHH$ is an embedding with $\varphi(S)=\Sigma$, then is $\varphi$ proper? i.e. whether the
preimage of a compact subset $K$ of $\BHH$, $\varphi^{-1}(K)$, is compact in $S$ for any $K$.

Anderson implicitly talks about this property in \cite{An1}, and \cite{An2}. Gabai conjectures the existence of a
properly embedded area minimizing plane in $\BH^3$ (and for any cocompact metric on $\BH^3$) for any given simple
closed curve $\Gamma$ in $\Si$. Later, Soma proved the existence of such an area minimizing plane in more general
situation (Gromov hyperbolic spaces) in \cite{So1} and \cite{So2}. Later, the author gave an alternative proof for
Soma's results in \cite{Co5}.

\begin{thm} \cite{So1}, \cite{So2}, \cite{Co5}
Let $X$ be a Gromov hyperbolic $3$-space with cocompact metric, and $\Si$ be the sphere at infinity
of $X$. Let $\Gamma$ be a given simple closed curve in $\Si$. Then, there exists a properly
embedded area minimizing plane $\Sigma$ in $X$ with $\PI \Sigma = \Gamma$.
\end{thm}

In recent years, the properly embeddedness of the complete minimal surfaces has also been in serious attack in $\BR^3$
case. This is called as Calabi-Yau Conjecture for minimal surfaces, and has been shown by Colding and Minicozzi in
\cite{CM2}. Later, the author showed an analogous result for hyperbolic space. In particular, he showed that for any
area minimizing plane $\Sigma$ in $\BH^3$ with asymptotic boundary $\Gamma$ which is a simple closed curve with one
smooth point, then $\Sigma$ is properly embedded in $\BH^3$. The technique is very different from Colding and
Minicozzi's techniques. While Colding-Minicozzi relates intrinsic distances and extrinsic distances for embedded
minimal surface in $\BR^3$ by using very powerful analytical techniques, the author's techniques are purely
topological.

\begin{thm} \cite{Co6}
Let $\Sigma$ be a complete embedded area minimizing plane in $\BH^3$ with $\PI\Sigma=\Gamma$ where
$\Gamma$ is a simple closed curve in $\Si$ with at least one smooth ($C^1$) point. Then, $\Sigma$
must be proper.
\end{thm}

\begin{pf} (Sketch) Assume that $\Sigma$ is a non-properly embedded area minimizing plane in $\BH^3$ with $\PI\Sigma=\Gamma$ where
$\Gamma$ is a simple closed curve in $\Si$ with at least one smooth point. The author gets a contradiction by analyzing
the disks in the intersection of $\Sigma$ with the balls $B_R(0)$ which exhaust $\BH^3$. First, he shows that for
sufficiently large generic $R>0$, $\Sigma\cap B_R(0)$ contains infinitely many disjoint disks. Then, he categorize
these disks as separating and nonseparating depending on their boundary in the annulus $A_R=CH(\Gamma)\cap \partial
B_R(0)$ is essential or not.

Then, he establishes the Key Lemma which shows that the nonseparating disks in $B_R(0)$ must stay close to the boundary
$\partial B_R(0)$. In particular, he proves that if $D_r$ is a nonseparating disk in $B_r(0)\cap\Sigma$, then there is
a function $F$ which is a monotone increasing function with $F(r)\rightarrow \infty$ as $r \rightarrow \infty$, such
that $d(0,D_r)>F(r)$ where $d$ is the distance. He proves the Key Lemma by using a barrier argument (See Figure 3). In
other words, by using the smooth point assumption, he proves the existence of a least area annulus $\mathcal{A}_r$ in
$\BH^3$ with $\PI \mathcal{A}_r = \Gamma_r^+\cup\Gamma_r^-$, where $\Gamma_r^\pm$ are simple closed curves sufficiently
close to $\Gamma$ in opposite sides. Since they are area minimizing, any nonseparating disk $D_r$ must stay in one side
of the least area annulus $\mathcal{A}_r$. As $r\to\infty$ the distance from $0$ to $\mathcal{A}_r$ will give the
desired function. Hence, this shows that nonseparating disks do not come close to $0$ point, and stay close to the
boundary $\partial B_r(0)$.

\begin{figure}[t]
\mbox{\vbox{\epsfbox{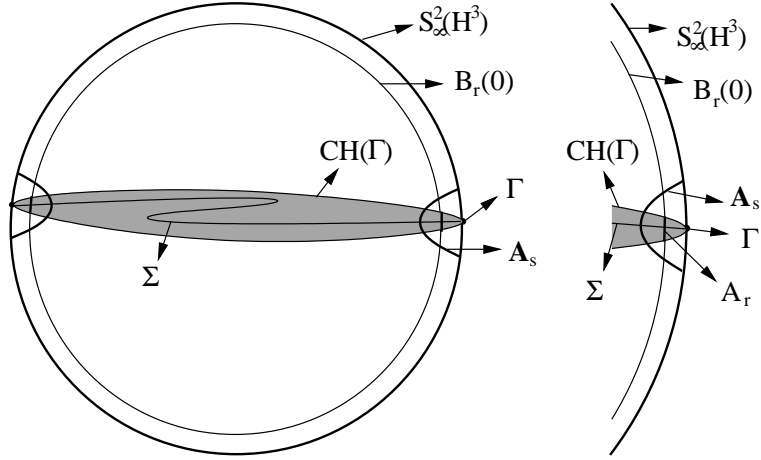}}} \caption{\label{fig:figure3} {The least area annulus $\mathcal{A}_s$ is used as
barrier in the proof of Key Lemma.}}
\end{figure}

Finally, the author proves the main result by using the Key Lemma as follows. A separating disk $D_{R_1}$ in
$\Sigma\cap B_{R_1}(0)$ will be a subdisk in a nonseparating disk $E_{R_2}$ in $\Sigma\cap B_{R_2}(0)$ where $R_2>R_1$.
By choosing $R_2$ appropriately and by using the fact that the separating disk $D_{R_1}$ is a subset of the disk
$E_{R_2}$, he shows that the nonseparating disk $E_{R_2}$ comes very close to $0$ point, which is a contradiction.
\end{pf}

We should add that it would be very interesting to generalize Colding-Minicozzi's result on Calabi-Yau conjecture for
minimal surfaces to this context \cite{CM2}, i.e. relating the intrinsic distances with the extrinsic distance for any
embedded minimal surface in $\BH^3$.

Note that properly embeddedness of absolutely area minimizing hypersurfaces in $\BHH$ is almost automatic. This is
because a nonproperly embedded area minimizing hypersurface in $\BHH$ would have an intersection of infinite volume
with a sufficiently large compact ball in $\BHH$, which is impossible for absolutely area minimizing hypersurfaces.
Also, in a forthcoming paper \cite{Co11}, the author constructs examples of non-properly embedded minimal planes in
$\BH^3$.

\subsection{The Global Structure}\ \\

On the other hand, the space of all solutions to the asymptotic Plateau problem is another interesting problem, and its
structure gives powerful global analysis tools to get important results on the number of solutions to the asymptotic
Plateau problem . In particular, the author showed that the space of minimal planes in $\BH^3$ with asymptotic boundary
$C^{3,\alpha}$ smooth simple closed curve is a manifold and its projection to the asymptotic boundary is a Fredholm map
in \cite{Co2}. By using these results, the author showed a generic uniqueness result (Theorem 5.4) for $C^3$ smooth
simple closed curves in $\Si$, \cite{Co3}.

Very recently, by using different techniques, Alexakis and Mazzeo generalized author's results to complete properly
embedded minimal surfaces of any fixed genus in convex cocompact hyperbolic $3$-manifolds ($\BH^3$ is a special case).

\begin{thm} \cite{AM}
Let $X$ be a convex cocompact hyperbolic $3$-manifold, and $M_k(X)$ is the space of properly embedded minimal surfaces
in $X$ of genus $k$ with asymptotic boundary $C^{3,\alpha}$ simple closed curve in $\PI X$. Let $\xi$ be the space of
all $C^{3,\alpha}$ curves in $\PI X$. Then, both $M_k(X)$ ($M_0(\BH^3)$ case in \cite{Co2}) and $\xi$ are Banach
manifolds, and the projection map $\pi_k: M_k(X) \to \xi$ is a smooth proper Fredholm map of index $0$.
\end{thm}

Note that being Fredholm map of index $0$ is a very strong property, and it can be considered as the map is locally
one-to-one for generic points. Indeed, they showed that $\pi_k$ is not only Fredholm of index $0$, but also proper.
Hence, by using this result, they developed a powerful $\BZ$-valued degree theory for $\pi_k$ as follows:

$$deg(\pi_k)=\sum_{\Sigma \in \pi_k^{-1}(\Gamma)} (-1)^{n(\Sigma)}$$

\noindent where $\Gamma$ is a regular value of $\pi_k$ and $n(\Sigma)$ is the number of negative eigenvalues of the
Jabobi operator $-L_\Sigma$. By combining this degree theory with the techniques in \cite{TT} and \cite{Wh2}, one can
get very interesting results on complete minimal surfaces in $\BH^3$ (see Section 4 in \cite{AM}).\\

\subsection{Intersections}\ \\

Another interesting property of the solutions to the asymptotic Plateau problem is that their
intersections mostly controlled by their asymptotic boundary. In many cases, if the asymptotic
boundaries are disjoint, then the solutions to the asymptotic Plateau problem are also disjoint.

\begin{thm} \cite{Co7}
Let $\Gamma_1$ and $\Gamma_2$ be two disjoint simple closed curves in $\Si$. If $\Sigma_1$ and
$\Sigma_2$ are area minimizing planes in $\BH^3$ with $\PI \Sigma_i = \Gamma_i$, then $\Sigma_1$
and $\Sigma_2$ are disjoint, too.
\end{thm}

The idea of the proof for this case is quite simple. If $\Sigma_1\cap\Sigma_2$ is not empty, then
as asymptotic boundaries are disjoint, the intersection must contain a simple closed curve
$\gamma$. Then, $\gamma$ bounds a disk $D_i$ in $\Sigma_i$. By swaping the disks, we get area
minimizing planes with a folding curve $\gamma$. Hence, we can reduce the area by smoothing out the
curve and get a contradiction \cite{MY2}. With slight modifications, this result can be generalized
to absolutely area minimizing hypersurfaces.

\begin{thm} \cite{Co7}
Let $\Gamma_1$ and $\Gamma_2$ be two disjoint connected closed codimension-$1$ submanifolds in
$\SI$. If $\Sigma_1$ and $\Sigma_2$ are absolutely area minimizing hypersurfaces in $\BHH$ with
$\PI \Sigma_i = \Gamma_i$, then $\Sigma_1$ and $\Sigma_2$ are disjoint, too.
\end{thm}

To generalize the idea of previous theorem, first by using the regularity result of Hardt-Lin in
\cite{HL} (Theorem 4.1), it can be showed that $\Sigma_1$ and $\Sigma_2$ are separating. As
asymptotic boundaries $\Gamma_1$ and $\Gamma_2$ are disjoint, the intersection $\alpha$ stays in
the compact part, and as $\Sigma_i$ separating, $\alpha$ separates a compact codimension-0 part
$S_i$ from $\Sigma_i$. Again by swaping these parts, one can get absolutely area minimizing
hypersurfaces with codimension-1 singularity set $\alpha$ which contradicts to interior regularity
results of geometric measure theory.

On the other hand, these arguments cannot be applied to minimal submanifolds, or area minimizing
submanifolds in a specified topological class. In the minimal submanifold case, the surgery
argument completely fails as there is no area factor to compare. The main problem with the fixed
topological class case is that the essential surgery argument in the proof is not working as after
surgery one may not stay in the same topological class. In the absolutely area minimizing case
there is no topological restriction.

There is a related conjecture which has important applications in $3$-manifold topology.

\vspace{0.3cm}

\noindent {\bf Disjoint Planes Conjecture:} Let $\Gamma_1, \Gamma_2$ be simple closed curves in
$S^2_\infty(X)$, where $X$ is a Gromov hyperbolic $3$-space with cocompact metric. If $\Gamma_1$
and $\Gamma_2$ do not cross each other (i.e. They are the boundaries of disjoint open regions in
$S^2_\infty(X)$), then any distinct area minimizing planes $\Sigma_1, \Sigma_2$ in $X$ with
asymptotic boundary $\Gamma_1, \Gamma_2$ are disjoint.

\vspace{0.3cm}

Even though this conjecture is interesting in its own right, it has powerful topological applications. The most
important application is constructing the area minimizing representative of an essential $2$-dimensional object in a
$3$-manifold, like incompressible surfaces, and genuine laminations. With this conjecture, if a $2$-dimensional
embedded essential object $S$ in a Gromov hyperbolic manifold $M$ induces a $\pi_1$-invariant family of circles $\PI
\widetilde{S}$ in $S^2_\infty(\widetilde{M})$, then by spanning the circles with area minimizing planes, the conjecture
would give you an $\pi_1$-invariant {\em pairwise disjoint} family of area minimizing planes in $\widetilde{M}$. Hence,
by projecting down the planes to $M$, it is possible to get the {\em embedded} area minimizing representative of $S$ in
the $3$-manifold $M$. Note that the author showed that the conjecture is generically true in \cite{Co9}.

\subsection{Renormalized Area}\ \\

In \cite{AM}, in addition to the study of the global structure of moduli spaces of complete minimal surfaces in $\BH^3$
and a $\BZ$-valued degree theory on them (see Section 7.2), Alexakis and Mazzeo defined a notion called {\em
renormalized area $\mathcal A(Y)$} for properly embedded minimal surfaces $Y$ in $\BH^3$ (or more generally convex
cocompact hyperbolic $3$-manifolds) where $\PI Y = \Gamma$ is a $C^{3,\alpha}$ simple closed curve in $\Si$. They
showed that if a minimal surface minimize renormalized area, it must be an area minimizing surface.

\begin{thm} \cite{AM} Let $\Gamma$ be a $C^{3,\alpha}$ embedded curve in $\Si$. Suppose that $Y_1$
and $Y_2$ are two properly embedded minimal surfaces in $\BH^3$ with $\PI Y_1 = \PI Y_2 = \Gamma$. If $Y_1$ is area
minimizing in $\BH^3$, then $\mathcal A(Y_1) \leq \mathcal A(Y_2)$, and equality holds if and only if $Y_2$ is also an
area minimizer.
\end{thm}

Moreover, they also showed that the renormalized area functional $\mathcal{A}$ is connected with the Willmore
functional $\mathcal{W}$, which is the total integral of the square of the mean curvature, in the following way. The
renormalized area functional is defined for any convex cocompact hyperbolic $3$-manifold $X$. After modifying the
metric on $X$ in a suitable way such that it induces a $Z_2$-invariant smooth metric on the double of $X$, say $2X$,
consider the double of any surface $\Sigma$ in $M_k(X)$ (see Section 7.2), say $2\Sigma$, in $2X$. Then, Alexakis and
Mazzeo showed that $\mathcal{A}(\Sigma)= -\frac{1}{2}\mathcal{W}(2\Sigma)$ for any $\Sigma\in M_k(X)$.

On the other hand, they also define an extended renormalized area $\mathcal R$ which is defined for all properly
embedded surfaces $Y$ which intersect $\Si$ orthogonally and $\PI Y = \Gamma$ is a $C^{3,\alpha}$ simple closed curve
in $\Si$. Then the extended renormalized area behaves just like the area for these surfaces.

\begin{thm} \cite{AM} Let $\Gamma$ be a $\mathcal C^{3,\alpha}$ closed curve in $\Si$.  Then the infimum of $\mathcal R(Y)$ where $Y$ ranges
over the set of all $\mathcal C^{3,\alpha}$ surfaces with $\PI Y = \Gamma$ which intersect $\Si$ orthogonally is
attained only by absolutely area-minimizing surfaces. Also, if $Y$ is a critical for $\mathcal R$, then $Y$ must be a
minimal surface.
\end{thm}

Notice that renormalized area behaves just like the area for these infinite surfaces in many ways. Hence, many
techniques from the compact area minimizing surfaces can be generalized to these surfaces with this new tool.

\vspace{0.3cm}

\begin{center}
{\bf Acknowledgements}
\end{center}

I would like to thank Urs Lang for very valuable remarks.

\end{document}